\documentclass[final]{beatcs}
\title{\bf Passages of Proof}
\author{{\bf Cristian S. Calude}$^{1}$,  \ {\bf Elena Calude}$^{2}$,\, {\bf Solomon Marcus}$^{3}$\\
$^{1}$University of
Auckland, New Zealand\\  {\tt  cristian@cs.auckland.ac.nz}\\
$^{2}$Massey University at Albany,
  New Zealand\\ {\tt e.calude@massey.ac.nz}\\
$^{3}$Romanian Academy, Mathematics, Bucharest, Romania\\  
{\tt Solomon.Marcus@imar.ro}}
\begin{document}
\date{\today}
\thispagestyle{empty}
\maketitle

\section{To Prove or Not to Prove--That  Is the Question!}
\hfill\mbox{\it   {\small Whether 'tis nobler in the mind to suffer}}\\[-3ex]

\hfill\mbox{\it   {\small The slings and arrows of outrageous fortune,}}\\[-3ex]

\hfill\mbox{\it   {\small Or to take arms against a sea of troubles}}\\[-3ex]

\hfill\mbox{\it    {\small And by opposing end them?}} \\[-3ex]

\hfill\mbox{Hamlet 3/1, by W. Shakespeare}
\bigskip

In this paper we propose a new  perspective on the evolution and history of the idea of mathematical
proof. Proofs will be studied at three levels: syntactical, semantical and pragmatical.
Computer-assisted proofs will be give a special attention.
Finally, in a highly speculative part,  we will anticipate the evolution of proofs under the assumption
that  the  quantum computer will materialize. We will argue that 
 there is little `intrinsic' difference between traditional and  `unconventional' types of proofs.

\section{Mathematical Proofs: An Evolution in Eight Stages}
\hfill\mbox{\it   {\small Theory is to practice as rigour is to vigour.}}   {\small D. E. Knuth}
\medskip

 {\it Reason} and   {\it experiment} are two ways to acquire knowledge.  For a long time
mathematical proofs  required only reason; this might be no longer true. 
   We can distinguish 
eight periods in the evolution of the idea of mathematical proof.
 The first period was that of  pre--Greek mathematics, for
instance the Babylonian one, dominated by observation,
intuition and experience.

 The second period was started by
Greeks such as Pythagoras and is characterized by the discovery
of deductive mathematics, based on theorems.  Pythagoras
proved his theorem, but  the respective statement was discovered
much earlier.  Deductive mathematics saw
 a culminating moment in Euclid's  geometry. The importance
of abstract reasoning to  ancient Greeks can be illustrated
by citing Aristophanes's comedy {\em The Birds}
which includes a cameo appearance of Meton, the astronomer, who
claims that he had squared the circle.  Knuth \cite{knuth} rhetorically
asked: ``Where else on earth would a playwright think of including such a scene\/?"
Examples would have been difficult to produce in 1985,
but today the situation has changed. Take for example, the movie {\em Pi}
written and  directed by Darren Aronofsky
Starring Sean Gullette or 
 Auburn's play {\em Proof}
\cite{auburn} originally produced by the Manhattan Theatre Club on 23rd
May 2000.

 In a more careful description, we observe that deductive mathematics
starts with Thales and Pythagoras, while the axiomatic approach begins with
Eudoxus and especially with Aristotle, who shows that a demonstrative
science should be based on some non--provable principles, some 
common to all sciences, others specific to some of them. Aristotle also used
 the expression  ``common notions"  for axioms (one of them being the famous
principle of non--contradiction). Deductive thinking and axiomatic thinking
are combined in Euclid's {\it Elements}  (who uses, like Aristotle, ``common
notions" for ``axioms"). The great novelty brought by Euclid is the fact
that, for the first time,  mathematical proofs (and, through them,
science in general)  are built on a long distance perspective, in a step
by step procedure, where you have to look permanently to  previous
steps and to fix your aim far away to the hypothetical subsequent steps.
Euclid became, for about two thousands years,  a term of reference for
the axiomatic--deductive thinking, being considered the highest standard
of rigour. Archimedes, in his treatise on static equilibrium, the
physicists of the Middle Age (such as Jordanus de Nemore, in 
{\it Liber de ratione ponderis}, in the 13th century), B. Spinoza in 
{\it Ethics}  (1677) and I. Newton in {\it Principia} (1687) follow Euclid's
pattern. This tradition is continued in many more recent works, not only
in the field of mathematics, but also in physics, computer science, biology, linguistics,
etc.
 
 However, some shortcomings of Euclid's approach were obstacles  for
the development of mathematical rigour. One of them was the fact that,
until Galilei, the mathematical language was essentially the ordinary
language, dominated by imprecision resulting from its predominantly
spontaneous use, where emotional factors and lack of care have an
impact. In order to diminish this imprecision and make the mathematical
language capable to face  the increasing need of precision and
rigour, the ordinary language had to be supplemented by an artificial
component of symbols, formulas and equations: with Galilei, Descartes,
Newton and Leibniz, the mathematical language became more and more a
mixed language, characterized by a balance between its natural and 
artificial components. In this way, it was possible to pack in a
convenient, heuristic way, previous concepts and results, and to refer 
 to them in the subsequent development
of mathematical inferences. To give only one example, one can imagine
how difficult was to express the $n$th power of a binomial expression
in the absence of a symbolic representation, i.e., using only words of the
ordinary language.  This was the third step in the development of
mathematical proofs.   
    
 The fourth step is associated with the so--called epsilon rigour, so
important in mathematical analysis; it occurred in the 19th century
and it is associated with names such as A. Cauchy and K. Weierstrass.
So, it became possible to renounce  the predominantly intuitive
approach via the infinitely small quantities of various orders, under the form
of functions converging in the limit to zero (not to be confused with the
Leibnizian infinitely small, elucidated in the second half of the 20th
century, by A. Robinson's non--standard analysis).
       The epsilon rigour brought by the fourth step  created
        the possibility to cope in a more accurate manner with
processes with infinitely many steps such as limit, continuity,
differentiability and integrability.

The %third
fifth
period begun with the end of the 19th century, when Aristotle's
logic, underlining   mathematical proofs for two thousands
years, entered a crisis with the challenge of the principle
of non--contradiction. 
This crisis
 was already announced by the discovery of
non--Euclidean geometries, in the middle  of the 19th century. 
Various therapies were proposed to
free the    mathematical proof of the dangerous effects of paradoxes (Russell--Whitehead,
Hilbert, Brouwer, etc). This period  covers the first
three decades of the 20th century and is dominated by the
optimistic view stating the possibility to arrange the whole
mathematics as a formal system and to decide  for any possible
statement whether it is true or false.  However, even during this period mathematicians were
divided with respect to the acceptance of non--effective (non--constructive) entities and proofs  (for example, 
Brouwer's intuitionism rejects the principle of excluded middle
in the case of infinite sets).  Intuitionism was a signal for
the further development of constructive mathematics, culminating
with the algorithmic approach leading to computer science.

 The  sixth period begins with G\"odel's incompleteness theorem
(1931), for many meaning  the unavoidable failure of any attempt
to formalise the whole of mathematics.  Aristotle's requirement
of complete lack of contradiction can be satisfied only by paying the price of
incompleteness of the working formal system. Chaitin (1975) has continued this
trend of results by proving that from $N$ bits of axioms one cannot
 prove that a program is the smallest possible if
it is more than $N$ bits long; he suggested
 that complexity is a source of
incompleteness  because a  formal system can capture only a tiny
amount of the huge information contained in the world of mathematical
 truth. This principle has been proved in Calude and J\" urgensen \cite{cj}.
 Hence,  incompleteness is natural and
inevitable   rather then mysterious and esoteric. This raises the natural question
(see Chaitin \cite{ch02}):
{\em How come that in spite of incompleteness, mathematicians are making so much
progress\/?}

 The   seventh period belongs to the second half of the 20th century,
when algorithmic proofs  become acceptable only when their
complexities were not too high. Constructiveness is no longer
enough,  a reasonable high complexity (cost) is mandatory.
We are now  living  in this period.  An important 
event of this period was   the  1976 proof of the Four--Colour
Problem (4CP): it marked
the reconciliation of empirical--experimental mathematics with
deductive mathematics, realized by the use of computer
programs as pieces of a mathematical proof. Computer refers
to classical von Neumann computer. At the horizon we can see
the (now hypothetical) quantum computer which may 
modify radically the relation between empirical--experimental mathematics and
deductive mathematics   \ldots

 With the eighth stage, proofs are no longer exclusively based on logic and
deduction, but also empirical and experimental. On the
other hand, in the light of the important changes brought by
authors like Hilbert, already at the beginning of the 20th
century, primitive terms became to have an explicit status, axioms
show their dependence on physical factors and the axiomatic--deductive
method displays its ludic dimension, being a play with abstract
symbols. Hilbert axiomatization of geometry is essentially different
from Euclid's geometry and this fact is well pointed out by 
 Dijkstra in \cite{dijkstra}  where he considers that, by directing their
attention towards provability, formalists 
 %took their distance in respect to
 circumvented  the vague metaphysical notion of ``truth".  Dijkstra
qualifies as  ``philosophical pollution" the mentality which pushed
 Gauss  not to  publish his ideas related to non--Euclidean
geometry. Contrary to  appearances, believes Dijkstra, Euclidean
geometry is not a prototype of a deductive system, because it is
based to a large extent on pictures (so--called definitions of points
and lines, for instance) motivated by the need of geometric
intuition. For Dijkstra, the claim 
 that the Euclidean geometry is a model of deductive thinking, is a
big lie. As a matter of fact, the shortcomings to which Dijkstra
refers were well-known, as  can be seen in Morris Kline's
 book \cite{kline}, pp. 86--88. 
In contemporary mathematics we are facing a change of perspective, a change of scenario, replacing the
old itinerary definition--theorem--proof by another one (see, for
instance, W. Thurston), based on ideas, examples and motivations.
The interesting fact is that the gap created between proof and
intuition by Hilbert prepared the way for a new marriage between
deduction and experiment, made possible by the computational
revolution, as it was shown by the latest step in the evolution of
proofs.

\section{Proofs, Theorems and Truths}
\hfill\mbox{\it   {\small Depuis les Grecs, qui dit Math\'ematique,}}\\[-3ex]

\hfill\mbox{\it    {\small dit d\'emonstration.}}   {\small Bourbaki}
\medskip

What is a mathematical proof\/? At a first glance the answer seems obvious: a proof  is a series of logical steps based on some axioms and deduction rules which  reaches a desired conclusion. Every step in a proof can be checked for correctness by examining it to ensure that it is logically sound. In David Hilbert's words:
``The rules
should be so clear, that if somebody gives you what they claim is a
proof, there is a mechanical procedure that will check whether the
proof is correct or not, whether it obeys the rules or not."
By making sure that every step is correct, 
 one  can tell once and for all whether   a theorem has been proved.  Simple! A moment of reflection shows that the problem may  not be so simple. For example,
what if the ``agent" (human or computer)  checking a proof for correctness makes a mistake (agents are fallible)\/? Obviously, another agent  has to check that the agent doing the checking did not make any mistakes.  Some other agent will need to check that agent, and so on. Eventually one runs out of agents who could check the proof and, in principle, they could all have made a mistake! 

 The mistake is the neighbour and the brother of proof, it is
both an opponent and a stimulus. An interesting analysis, responding to Joseph L. Doob's challenge, of various possible mistakes
in the proof of the 4CT can be found in the
 work of  Schmidt \cite{schmidt}.
In 1976, Kenneth Appel and Wolfgang Haken  proved the 4CT. They used some of Alfred  Kempe's ideas, but avoided his mistake.\footnote{In 1879 Kempe announced his `proof' of the 4CT in both the magazine {\em Nature} and the {\em American Journal of Mathematics}. Eleven years later, Percy Heawood found an error in the proof which nobody had spotted, despite careful checking.} They showed that if there is  a map which needs five colours, then a contradiction
follows. If there are several five--colour maps, they have chosen one with the smallest number of countries and proved that this map must contain one of 1,936 possible configurations; they also proved that every one of these possible configurations can be reduced into a smaller configuration which also needs five colours. This is a contradiction because we assumed that we already started with the smallest five--colour map. The reduction step, i.e., the step in which one shows that the 1,936 configurations could be reduced was actually done  by brute force computer search through every configuration. No human being could  ever actually read the entire proof to check  its correctness. For Ron Graham, ``The real question is this: If no human being can ever hope to check a proof, is it really a proof\/?" 

In 1996  Robertson,  Sanders, Seymour and  Thomas \cite{rsst} offered a simpler proof involving only 633 configurations.  The paper \cite{rsst} concludes with the following interesting comment (p. 24):
 ``We should mention  that both our programs use only integer
arithmetic, and so we need not be concerned with round--off errors
and similar dangers of floating point arithmetic.
However, an argument can be made that our ``proof" is not a proof in
the traditional sense, because it contains steps that can never be
verified by humans. In particular, we have not proved the 
correctness of the compiler
we compiled our programs on, nor have we proved the infallibility of the
hardware we ran our programs on. These  
have to be taken on faith, and are conceivably a source of error.
However, from a practical point of view, 
the chance of a computer error that appears consistently
in exactly the same way on all runs of our programs on all the compilers
under all the operating systems that  our programs run on is 
infinitesimally small compared to the chance of a human
error during the same amount of case--checking. Apart from
this hypothetical possibility of a computer consistently giving
an incorrect answer, the rest of our proof can be verified in the
same way as traditional mathematical proofs. We concede, however, that
{\it verifying a computer program is much more difficult than checking
a mathematical proof of the same length}."\footnote{Our Italics.}

According to Vladimir Arnold, ``Proofs are to  mathematics what spelling (or even calligraphy) is to 
poetry. Mathematical works do consist of proofs, just as 
 poems do consist of characters." 
These analogies  point  out both  the necessity and the
insufficiency of proofs in the development of mathematics.  Indeed,
spelling is the way poetry takes expression, but it is equally
the tool used by the common everyday language, in most cases
devoid of any poetic effect. What should be added to spelling
in order to get a piece of poetry remains a mystery. A poem
consists of characters, but it is much more than a meaningful
concatenation of characters.

Mathematics cannot be conceived in the absence of proofs.
According to  Foia\c s \cite{foias}, ``the theorem is the brick of
mathematics". Obviously, ``proof" and  ``theorem" go together;
the object of a proof is to reach a theorem, while theorems are
validated by proofs. Theorems are, for the construction of
mathematics, what bricks are   for the construction of a
building. A building is an articulation of bricks and,
analogically,  a mathematical work is an articulation of theorems.
Motivated by a similar view, Jean Dieudonn\'{e} \cite{dieu} defines a
mathematician as a person who has proved at least one theorem.
In contrast, Arnold's analogies
point out the fact that mathematics is much more than a chain
of theorems  and proofs, so implicitly a mathematician should be
much more than the author of a theorem. Probably the best example
is offered by Bernhard Riemann whose lasting fame does not come (in the first instance)
from his theorems or proofs, but from his conjectures, definitions, concepts
and examples (see for example, the discussion in 
Hersh \cite{hersh}, pp. 50--51). Srinivasa Ramanujan is another famous example of a  mathematician
who produced more results than proofs. What the mathematical community seems to value
most are ``ideas". ``The most respected mathematicians are those with strong
`intuition' " (Harris \cite{harris}, p. 19).

\section{Mathematical Proofs: The Syntactic Dimension}
\hfill\mbox{\it %Rigour  but surely not rigour mortis.
  {\small Cum Deus calculat,   fit mundus.}}   {\small Leibniz}
\medskip

Of course, the first thing to be discussed is G\"{o}del's incompleteness theorem (GIT)
which says that {\it every formal system which is (1) finitely specified, (2) rich enough
to include the arithmetic, and (3) consistent, is incomplete.}  That is, there exists
an arithmetical statement which  (A) can be expressed in the formal system, (B) is  true, but 
(C) is unprovable within the formal system. All conditions are necessary.
 Condition (1) says that there is an
algorithm listing all axioms and inference rules (which could be infinite). Taking
as axioms all true arithmetical statements will not do, as this set is not finitely
listable. But what does it mean to be a ``true arithmetical statement"? It is a statement
about non-negative integers which cannot be invalidated by finding any combination
of non-negative integers that contradicts it. In Alain Connes terminology (see
\cite{cls}, p. 6), a true arithmetical statement is a ``primordial mathematical reality".
Condition (2) says that the formal system has all the symbols and axioms used in arithmetic,
the symbols for  $0$ (zero), $S$ (successor), $+$ (plus),  $\times$ (times), $=$ (equality)
and the axioms making them work (as for example, $x +S(y) = S(x+y)$). Condition (2)
cannot be satisfied if you do not have individual
terms for $0, 1, 2, \dots $; for example, Tarski 
proved that  Euclidean geometry, which refers to points, circles and lines,
 is complete. Finally (3) means that the formal system is free of contradictions. The essence of GIT is to distinguish
between truth and provability. A closer  real life analogy  is the distinction between
truths and judicial decisions, between what is true and what can be proved in court.\footnote{The Scottish
judicial system which admits three forms of verdicts, guilty, not--guilty and not--proven, comes
closer to the picture described by GIT.}  How large is the set of true
and unprovable statements? If we fix a formal system satisfying
all three conditions in GIT, then the set of true
and unprovable statements is topologically ``large" (constructively, 
a set of second Baire category, and in some cases even ``larger"), cf.
Calude, J\" urgensen, Zimand \cite{cjz}; because  theorems proven in such a system  have bounded complexity,   the probability that  an  $n$-bit statement is provable tends to zero when $n$ tends to infinity (see Calude and J\" urgensen \cite{cj}).

 There is a variety of reactions in interpreting GIT,
ranging from pessimism to optimism or simple dismissal (as irrelevant for the practice of
mathematics). For pessimists, this result can be interpreted as
the final, definite failure of any attempt
to formalise the whole of mathematics. For example, 
Hermann Weyl acknowledged that  GIT
has exercised a ``constant drain on the enthusiasm" with which he has engaged 
himself in mathematics and for Stanley Jaki, GIT is a fundamental
barrier in understanding the Universe. In contrast,  scientists like
Freeman Dyson acknowledge the limit placed by GIT  on our
ability to discover the truth in mathematics, but interpret this in an
optimistic way, as a guarantee that mathematics will go on forever (see  
Barrow \cite{barrow}, pp. 218--221).

In  modern times a penetrating insight into the incompleteness phenomenon
has been obtained by an information--theoretic analysis pioneered by Chaitin
in \cite{ch75}.   Striking results have been obtained by studying
the Chaitin's Omega Number, $\Omega$, the halting probability of a self-delimiting universal 
Turing machine. This number is not only uncomputable, but also (algorithmically)
random. Chaitin has proven the following important theorem:
 {\it If $ZFC$  (Zermelo set theory with the Axiom of Choice)  is arithmetically
sound,  that is, any theorem of arithmetic proved by $ZFC$ is \emph{true}, 
then,  $ZFC$ can determine
the value of only finitely many bits of $\Omega$, and one can
give a bound on the number of bits of
$\Omega$ which
$ZFC$ can determine.} 
 Robert Solovay \cite{solovay2k} (see more in \cite{cc,crisomega,cris,cris2002})
 has  constructed {\it a self-delimiting universal 
Turing machine
 such that
$ZFC$, if arithmetically sound, cannot determine any single bit of  its halting probability}
($\Omega$).
Re--phrased, the most powerful formal axiomatic system  is powerless when dealing
with the questions of the form 
 ``is the $m$th bit of  $\Omega$  0\/?" 
or 
``is the $m$th bit of $\Omega$  1\/?". 

 Chaitin has constructed an exponential Diophantine equation
$F(t; x_1,   \ldots   ,x_n)=0$ with the following property: the infinite
binary sequence whose $m$th term is 0  or 1 depending whether the
equation $F(m; x_1,   \ldots   ,x_n)=0$ has finitely or infinitely many
solutions is  exactly the digits of $\Omega$, hence it is
random; its infinite
amount of information  is algorithmically
incompressible. The importance of exponential Diophantine equations comes from the fact that
most problems in mathematics can be formulated in terms of   these type of equations;
Riemann's Conjecture is one such example.
 Manin \cite{manin1}, p. 158, noticed that ``The epistemologically
important point is the discovery that randomness can be defined without
any recourse to physical reality   \ldots  \, in such a way that the necessity
to make an infinite search to solve a parametric series of problems leads
to the technically random answers. Some people find it difficult to
imagine that a rigidly determined discipline like elementary arithmetic
may produce such phenomena".

Last but not least, is the truth achieved through a formal proof the ultimate
expression of knowledge\/? Many (mathematicians) will give a positive
answer, but perhaps not all. For the  13th century Oxford
 philosopher Roger Bacon,  ``Argument reaches a conclusion and compels us to admit it, but it neither makes us certain nor so it annihilates
doubt that the mind rests calm in the intuition of truth, unless
it finds this certitude by way of experience." More recently, I. J. 
Schoenberg\footnote{Landau's son-in-law.}
is cited by Epstein (\cite{hahn}) as  saying that Edmund Landau
kept in his desk drawer for years a manuscript proving what is now called
 the two constants theorem: he had the complete proof  but could not 
believe it until his intuition was ready to accept it. Then he published it.
A ``proof is only one step in the direction of confidence" argued De Millo,
Lipton and Perlis in a classical paper on proofs, theorems and programs \cite{demillo}.
Written in the same spirit is
Don Knuth's warning: ``Beware of bugs in the above code: I have only proved it
correct, not tried it." 

\section{Mathematical Proofs: The Semantic Dimension}
\hfill\mbox{\it   {\small If one must choose between rigour and meaning,}} \\[-3ex]

\hfill\mbox{\it   {\small I shall unhesitatingly choose the latter.}}   {\small R. Thom}
\medskip

The above quotation turned slogan as ``more rigour, less
meaning", or better still, ``less rigour, more meaning"
(Chaitin  \cite{gregpccris})
points out the necessity to distinguish between the syntactic
and the semantic aspects of proofs.  Should proofs belong
exclusively to logic, according to the tradition started
by Greeks such as Pythagoras and Euclid\/? Or should they also be
accepted  as a cocktail of logical and empirical--experimental
arguments, as in the proof of the 4CT (1976)\/?  Mathematicians are
now divided into those giving an affirmative answer to the first question 
and implicitly a negative answer to the
second question and  those giving a
negative answer to the first question and  an affirmative  one to the second question.
Computationally oriented mathematicians  usually belong to the
second category, while many other mathematicians (as, for
instance, the Fields medalist William Thurston) belong to  the
first, so  for them, the 4CT is not yet
proved! 
Meaning is a key distinction. For mathematicians such as Ren\'{e} Thom, Daniel Cohen and William
Thurston, correctness by itself does not validate a proof; 
it is also necessary to ``understand"
it. ``The mission of mathematics is understanding" says Cohen.
  Paul Halmos has also
insisted on the ``conceptual understanding".
 For him a ``good" proof of a theorem is one that sheds light on why it is true. It is just the process of understanding which is in question
with proofs like that given to the 4CT. Referring to
the proof of the 4CT, Halmos says: ``I do not find it easy to say what we learned from all that.
  \ldots \,   The present proof relies in effect on an Oracle, and I say down with Oracles!
They are not mathematics!"  In contrast with Halmos, who hopes that ``100 years from now
the map theorem will be   \ldots \,  an exercise in a first--year graduate course, provable in a couple of pages 
by means of appropriate concepts, which will be completely familiar by then" (see \cite{hersh}, p. 54),
R.~Hersh  thought that the problem itself might be responsible for the way it was solved: he is cited  by saying dejectedly
``So it just goes to show, it wasn't a good problem after all" (see \cite{casti} p. 73).

We will return later to these issues. For the moment
we make the following two observations.

\begin{description}
  \item[ A)] Not only the hybrid
proofs obtained as a combination of logical and empirical--experimental
arguments might be hard/impossible to be understood in their ``globality"; this happens   also for some pure deductive proofs.
An example is the proof  of the classification of finite simple groups called by Danny Gorenstein 
 the ``Thirty
Years War" (for the classification
battles were fought mostly in the
decades 1950--1980),
a work which comprises about 10,000--15,000 pages scattered in 500 journal articles by some 100 authors.\footnote{Still, there is a controversy in the mathematical community on whether these articles provide a complete and correct proof. For a recent account see Aschbacher \cite{ma}.}

  According to Knuth \cite{knuth} p. 18, ``\ldots \, program--writing is substantially more demanding than book--writing".  ``Why is this so\/? I think the main reason is that a larger attention span is needed when working on a large computer program than when doing other intellectual tasks.
\ldots \, Another reason is \ldots \, that programming demands a significantly higher standard of accuracy. Things don't simply have to make sense to another human being, they must make sense to a computer." Knuth compares his  \TeX  \, compiler (a document of about 500 pages)
with Feit and Thompson \cite{ft} theorem that all simple groups of odd order are cyclic. He lucidly argues that the program might not incorporate as much creativity and ``daring" as the proof of the theorem, but they come even when compared on  depth of details, length and paradigms involved. What distinguishes the program from the proof is the ``verification": convincing a couple of (human) experts that the proof  {\it works in principle}
seems to be easier than making sure that the program  {\it really works}. A demonstration that {\it there exists a way to compile \TeX \,}  is not enough! 
Another example, which will be discussed later in this section, is the proof of   Fermat's Last Theorem (FLT).  

\item[B)]  Without diminishing in any way the  ``understanding" component of mathematics we 
note that the idea of distinguishing between  ``good" and ``bad"  proofs on the light they shed on their
own truth seems to be, at least to some extent, relative and subjective.
\end{description}

%\section{Proofs Under the Status of Conjugate Pairs}

Thom's slogan  `more rigour, less meaning' was the main
point in his controversy with Jean Dieudonn\'e (as a
representative of the Bourbaki group). Taking rigour as
something that can be acquired only at the expense of
meaning and conversely, taking meaning as something
that can be obtained only at the expense of rigour,
we oblige mathematical proof to have the status of
what was called in physics a  ``conjugate (complimentary) pair", i.e.,
a couple of requirements, each of them being satisfied
only at the expense of the other (see \cite{marcus}).
 Famous prototypes of conjugate pairs are (position,
momentum) discovered by W. Heisenberg in quantum mechanics
and  (consistency, completeness)  discovered by K. G\" odel in
logic.
 But similar warnings come from other directions.
According to Einstein (see, for instance, \cite{rosen} p. 195),
  ``in so far as
the propositions of mathematics are certain, they do not
refer to reality, and in so far as they refer to reality,
they are not certain", hence (certainty,
reality) is a conjugate pair. Obviously, reality is here understood
as an empirical entity,    hence  mixed with all kinds
of imprecision, ranging from obscurity and disorder to
randomness, ambiguity and fuzziness \cite{marcus1}. 
   Pythagoras' theorem is
certain, but its most empirical tests will fail. There are
some genuine obstacles in our attempts to eliminate or at
least to diminish the action of various sources of imprecision.
Einstein implicitly calls our attention on one of them.  Proof,
to the extent to which it wants to be rigorous, to give us
the feeling of certainty, should be mathematical; but satisfying
this condition, means failing to reach reality. In other words,
the price we have to pay to obtain proofs giving us
the feeling of total confidence is to renounce to be directly
connected to reality. There is a genuine tension between
certainty and reality, they form a conjugate pair, which is
the equivalent of what  in the field of humanities is an
oxymoronic pair. However, there is an essential difference
between G\"odel's conjugate pair (consistency, completeness)
and Einstein's conjugate pair (certainty, reality). While
consistency and completeness are binary logical predicates,
certainty and reality are a matter of degree, exactly like
the terms occurring in Thom's conjugate pair: rigour and meaning.
In last two situations there is
room for manipulation and compromise.
   
 Near to the above conjugate pairs is a third one: (rigour,
reality), attributed to Socrates (see \cite{renyi}). 
   A price we have to pay   in order to reach
rigour is the replacement of the real world by a fictional 
one. There is no point and no line in the real world, if we
take them according to their definitions in Euclid's {\em Elements}.
Such entities belong to a fictional/virtual universe, in the same way
in which the characters of a theatrical play are purely
conventional, they don't exist as real persons. The rules of
deduction used in a mathematical proof belong to a game in
the style they are described in the scenario of a Hilbert
formal system, which is, as a matter of fact, a machine
producing demonstrative texts. A convention underlines the
production of theorems and again a convention is accepted 
in a theatrical play. In the first case, the acceptance of the
convention is required from both the author of the proof and
its readers;  in the second case all people involved, the author,
the actors and  spectators, have to agree  the proposed convention.
 Since many proofs, if not most of them, are components of a
modeling process, we have to add the unavoidable error of
approximation involved in any cognitive model. The model should satisfy
opposite requirements, to be as near as possible to the phenomenon
modelled, in order to be relevant; to be as far as possible from the
respective phenomenon, in order to useful, to make possible the existence of
at least one method or tool that can be applied to the model,
but not to the original (see \cite{marcus2}).
Theorems are discovered, models are invented. Their interaction
leads to many problems of adequacy, relevance and correctness,
i.e., of syntactic, semantic and pragmatic nature.
      
 In the light of the situations pointed out above, we can understand
some ironical comments about what a mathematician could be. It is
somebody who can prove theorems, as Dieudonn\'e claimed. But what kind
of problems are solved in this way\/?  ``Any problem you want, \ldots  except
those you need", said an engineer, disenchanted by his collaboration with
a mathematician. Again, what is a mathematician\/?  ``It is a guy capable
to give, after a long reflection, a precise, but useless answer",  said
another mathematician with a deep feeling of self irony. Remember the
famous reflection by Goethe:  ``Mathematicians are like French people,  they
take your question, they translate it in their language and you no longer
recognize it".

But things are controversial even when they concern   syntactic 
correctness. In this respect, we should distinguish two
types of syntactic mistakes: benign and malign. Benign
mistakes  have only a   local, not
global effect: they can be always corrected. Malign mistakes,
on the contrary, contaminate the whole approach and invalidate
the claim formulated by the theorem. 
When various authors
(including the famous probabilist J. L. Doob, see \cite{schmidt}) found some mistakes
in the proof of the 4CT, the authors of the
proof succeeded in showing that all of them were benign and  more
than this, {\it any other possible  mistake, not yet discovered,
should be benign}. How can we accept such arguments, when
 the process of  global understanding  of the respective
proof is in question\/? The problem remains open. 
A convenient, but fragile, solution is to accept  Thom's
pragmatic proposal: a theorem is validated
if it has been accepted by a general agreement\footnote{Perhaps  ``general" should
be replaced here by ``quasi--general".} 
of the mathematical community (see \cite{thom1,thom2}).

 The
problems raised by the 4CT were discussed by
many authors, starting with Tymoczko \cite{tymoczko} and    Swart \cite{swart} (more recent
publications are D. MacKenzie \cite{mc},
J. Casti \cite{casti}, A.S. Calude \cite{andreea}).
 Swart proposed the
introduction of a new entity called {\em agnogram}, which is ``a
theorem--like statement that we have verified as best we could,
but whose truth is not known with the kind of assurance we
attach to theorems and about which we must thus remain, to
some extent, agnostic." There is however the risk to give
the status of agnogram to any property depending on a
natural number $n$ and verified only for a large, but finite
number of values of $n$. This fact would be in conflict with 
Swart's desire  to consider an agnogram less than a
theorem, but more than a conjecture. Obviously, the 4CT is for Swart an agnogram, 
not a theorem. What is missing from
an agnogram to be a theorem\/? A theorem is a statement which could be 
checked individually by a mathematician and confirmed also
individually by at least two or three more mathematicians,
each of them working independently. But already here we can
observe the weakness of the criterion: how many
 mathematicians  are to check individually and
independently the status of an agnogram to give it the status of theorem\/?

The seriousness of this
objection can be appreciated by examining the case of Andrew Wiles' proof of FLT---a challenge to mathematics  since 1637 when Pierre de Fermat wrote it into the margin of one of his books.
 The proof is extremely intricate,  quite long (over 100 printed  pages\footnote{Probabilists would argue that very long proofs can at best be viewed as only probably correct, cf. \cite{demillo}, p. 273.  In view of \cite{cj}, the longer the statement, the lesser  its chance is to be proved.}), and  only a handful of people in the entire world can claim to  understand it.\footnote{Harris \cite{harris} believes that no more than 5\% of mathematicians  have made the effort to work through the proof. Does this have anything to do with what
  George Hardy has noted in   his famous {\em Apology}: ``All physicists
and a good many quite respectable mathematicians are contemptuous about proof."\/?} To the rest of us, it is utterly incomprehensible, and yet we all feel entitled  to  say  that ``the FLT has been proved". On which grounds\/? We say so because {\it we believe the experts} and {\em we cannot tell for ourselves}. Let us also note that in the first instance the original 1993 proof seemed accepted, then a gap was found, and finally it took Wiles and Richard Taylor another year to fix the error.\footnote{According to Wiles, ``It was an error in a crucial part of the argument, but it was something so subtle that I'd missed it completely until that point. The error is so abstract that it can't really be described in simple terms. Even explaining it to a mathematician would require the mathematician to spend two or three months studying that part of the manuscript in great detail."}

According to Hunt \cite{hunt}, ``In no other field of science would this be good enough. If a physicist told us that light rays are bent by gravity, as Einstein did, then we would insist on experiments to back up the theory. If some biologists told us that all living creatures contain DNA in their cells, as Watson and Crick did in 1953, we wouldn't believe them until lots of other biologists after looking into the idea agreed with them and did experiments to back it up. And if a modern biologist were to tell us that it were definitely possible to clone people, we won't really believe them until we saw solid evidence in the form of a cloned human being. Mathematics occupies a special place, where we believe anyone who claims to have proved a theorem on the say---so of just a few people---that is, until we hear otherwise."

Suppose we loosely define a religion as any discipline whose foundations rest on an element of faith,
irrespective of any element of reason which may be present. Quantum mechanics, for example, would qualify as a
religion under this definition. Mathematics would hold the unique position of being a branch
of theology possessing a ``proof" of the fact that it should be so classified.
``Where else do you have absolute truth\/? You have it in mathematics and you have it in religion, at least for some people. But in mathematics you can really argue  that this is as  close
to absolute truth as you can get" says Joel Spencer.
\section{%The Role of Proofs in Modelling
Mathematical Proofs: The Pragmatic Dimension}
\hfill\mbox{\it   {\small  Truth is not where you find it, but where you put it.} }   {\small  A. Perlis}
\medskip

 In the second half of the 20th century, theorems together with their
proofs occur with increasing frequency as components of some
cognitive models, in various areas of knowledge. In such
situations we are obliged to question the theorems not only
with respect to their truth value, but also in respect to their
adequacy and relevance within the framework of the models to
which they belong. We have to evaluate the explanatory capacity
of a theorem belonging to a model B, concerning the phenomenon
A, to which B is referring. This is a very delicate and
controversial matter, because adequacy, relevance and
explanatory capacity are a matter of degree and quality,
which cannot be settled by binary predicates. Moreover, there
is no possibility of optimization of a cognitive model. Any
model can be improved, no model is the best possible. This
happens because, as we have explained before,
 a cognitive model B of an entity A  has simultaneously  the tendency to increase its similarity with A
and stress its difference from A. To give only one example in this respect, we recall
the famous result obtained by  Chomsky \cite{chomsky}, in the late 1950s, 
 stating that context--free grammars are not able to generate the
English language. This result was accepted by the linguistic
and computer science communities until the eighties, when new
arguments pointed out the weakness of Chomsky's argument; but
this weakness was not of a logical nature, it was a weakness in
the way we consider the entity called ``natural language". As a
matter of fact, the statement ``English is a context--free
language" is still controversial.

Mathematical proofs are ``theoretical" and ``practical". Theoretical
proofs (formal, ideal, rigorous) are models for practical proofs (which are informal, imprecise,
incomplete). ``Logicians don't tell mathematicians what to do. They make
a theory out of what mathematicians actually do", says  Hersh
\cite{hersh}, p. 50. According to the same author, logicians study
what mathematicians  do the way fluid dynamicists study water waves.
Fluid dynamicists don't tell water how to wave, so logicians don't tell 
mathematicians what to do. The situation is not as simple as it appears.
Logical restrictions and formal models (of proof)  can play an important role 
in the  practice of mathematics. For example, 
the key feature of constructive mathematics is the identification
``existence =  computability" (cf. Bridges \cite{bridges}) and
a whole variety of constructive mathematics, the so--called Bishop
 constructive mathematics, is   mathematics with
intuitionistic rather than  classical underlying logic.

\section{Quasi--Empirical Proofs: From Classical to Quantum}
\hfill\mbox{\it   {\small
Truth does not change because it is, or is not, believed
}} \\[-3ex]

\hfill\mbox{\it   {\small by a majority of the people.}}   {\small Giordano Bruno}
\medskip

The use of large--scale programs, such as Mathematica,  Maple or  MathLab
is now widespread for symbolical and numerical calculations as well as  for 
graphics and  simulations.  To get a feeling of the extraordinary power of  such programs 
 one can visit, for example,  the Mathematica website {\tt http://www.wolfram.com}. 
New other  systems are produced; ``proofs as programs", ``proof animation"
or ``proof engineering" are just a few examples (see \cite{hayashi}). 
In some cases an experiment  conveys an aesthetic appreciation of
mathematics  appealing to a much broader audience (cf. \cite{bb1,bbg,crismarcus}).
 A significant, but simple example of the role an experiment
may play in a proof is given by  Beyer \cite{beyer}. He refers 
to J. North who asked 
for a computer demonstration 
that the harmonic series diverges. We quote Beyer: ``His example
illustrates the following principle: Suppose that one has a
computer algorithm alleged to provide an approximation to some
mathematical quantity. Then the algorithm should be accompanied
by a theorem giving a measure of the distance between the output
of the algorithm and the mathematical quantity being approximated.
For the harmonic series, one would soon find that the sum was
infinite." It is interesting to mention that  in 1973 Beyer made together
with Mike Waterman a similar attempt to compute Euler's constant;
 their experiment failed, but the error was discovered later by
 Brent \cite{brent}.

New types of proofs motivated by  the  experimental  ``ideology'' have appeared. For example, rather than being a static object, the {\it interactive proof} (see  Goldwasser, Micali, Rackoff
\cite{GMR}, Blum \cite{Blum})  is a two--party protocol  in which the {\it prover} tries to prove a certain fact to the {\it verifier}. During the interactive proof  the {\it  prover} and the {\it  verifier} exchange messages and at the end the {\it verifier} produces a verdict ``accept" or ``reject". 
A holographic (or probabilistic checkable) proof (see Babai \cite{Babai})  is still a static object but it is verified probabilistically. Errors become almost instantly apparent after a small part of the proof was checked.\footnote{More precisely, a traditional proof of length $l$ is checked in time a  constant power of $l$ while a holographic proof requires  only constant power of $\log_{2}l$. To appreciate the difference, the binary logarithm of the number of atoms in the known Universe is smaller than 300.} The transformation of a classical proof
(which has to be self-contained and  formal) into a holographic one requires super-linear time. 

The blend of logical and empirical--experimental
arguments (``quasi--empirical mathematics" for 
 Tymoczko \cite{tymoczko},
Chaitin \cite{ch00,ch02, gregphil}  or ``experimental mathematics" for Bailey,  Borwein \cite{bb}, Borwein, Bailey  \cite{bb1}, Borwein,  Bailey,   Girgensohn \cite{bbg}) may
lead to a new way to understand (and practice) mathematics. For example,
  Chaitin argued that
we should introduce the Riemann hypothesis
as an axiom: ``I believe that elementary number theory and the rest of mathematics should
be pursued more in the spirit of experimental
science,  and  that  you  should  be  willing  to
adopt new principles. I believe that Euclid's
statement that an axiom is a self--evident truth
is a big mistake. The Schr\"{o}dinger equation
certainly isn't a self--evident truth! And the Riemann hypothesis 
isn't self--evident either, but
it's very useful. A physicist would say that there
is ample experimental evidence for the Riemann hypothesis 
and would go ahead and take
it as a working assumption."
Classically, there are two  equivalent ways to look at the mathematical
notion of proof:  {\it logical}, as a finite sequence of sentences strictly obeying some
axioms and inference rules, 
and  {\it computational},  as a specific type of computation. 
Indeed, from a proof given as a sequence
of sentences one can easily construct a Turing machine producing that sequence as the
result of some finite computation and, conversely,
 given a machine computing
a proof we can just print all sentences produced during the computation and
arrange them into a sequence.

  This gives mathematics an immense advantage over
any science: a proof is an explicit sequence of reasoning steps that
can be inspected at {\it leisure}.  {\it  In theory}, if followed with care, such a sequence
either reveals a gap or mistake, or can convince a sceptic of its conclusion, 
in which case the theorem {\it is considered proven}. 
The equivalence between the logical and computational proofs has stimulated the
construction of programs which play the role of  {\it ``artificial" mathematicians}. 
 The
``theorem provers" have been very successful as ``helpers"  in proving many results, from
simple theorems of Euclidean geometry to the computation of a few digits of a Chaitin
Omega Number \cite{crisds}.
 ``Artificial" mathematicians are far less ingenious and subtle than human 
mathematicians, but 
they surpass their human counterparts by being infinitely more patient and diligent.

If a conventional proof is replaced by an ``unconventional" one (that is a proof consisting of
a sequence of reasoning steps   obeying
axioms and inference rules which depend not only on some logic, but also on the external
physical medium),  then the
conversion from a computation to a sequence of sentences may be impossible,
e.g. due to the size of the computation.
An extreme, and for the time being hypothetical  example, is the proof obtained 
as a result of a quantum computation (see Calude and P\u{a}un \cite{cp}). 
The quantum automaton would say ``your conjecture is true", but (due to  quantum interference) there
will be no way to exhibit all trajectories followed by the quantum automaton
in reaching that conclusion. 
The quantum automaton has the ability
to check a proof, but it may fail to reveal any ``trace" of the proof for the
human being operating the quantum automaton. 
Even worse, any attempt to
{\it watch} the inner working of the quantum automaton (e.g. by ``looking" inside
at  any information concerning the state of the ongoing proof) may
compromise forever the proof itself!
We seem to go back  to  Bertrand Russell
who said that %Thus 
``mathematics may be defined as the subject in which we 
                  never know what we are talking about, nor whether what we are 
                  saying is true", and even beyond by adding
{\it  and even when it's true we might not know why.}

Speculations about quantum proofs  {\it may not affect} the essence of mathematical objects and constructions
(which, many believe,  have an autonomous reality quite independent of the physical reality), but they
seem to {\it have an impact} on how we
 {\it learn/understand mathematics,} 
which is through
the physical world. 
Indeed, our glimpses of mathematics are revealed only
through physical objects, human brains, silicon computers, quantum automata, etc.,
hence, according to  Deutsch \cite{deutsch-97}, they have to obey not only the axioms and the
 inference  rules of the theory, but the {\it laws of physics} as well. To complete the picture we need to take into account also the {\it biological} dimension. No matter how precise the rules
 (logical and physical) are,  we need human consciousness to apply the rules and to understand
 them and their consequences. Mathematics is a human activity.

 \section{Knowledge Versus Proof}

 \hfill\mbox{\it 
  {\small The object of mathematical rigour is to sanction and}} \\[-3ex]

\hfill\mbox{\it  {\small  legitimize the conquests of intuition.}}  {\small J. Hadamard}
\medskip

 Are there intrinsic differences between traditional and `unconventional'
 types of proofs? 
  To answer this question we will consider the following
 interrelated questions:\\[-4ex]
 \begin{enumerate}
 \item Do `unconventional' methods supply us with a proof in some formal language?\\[-4ex]
 \item Do `unconventional' methods supply us with a mathematical proof?\\[-4ex]
 \item Do `unconventional' methods supply us with knowledge?\\[-4ex]
 \item Does mathematics require knowledge or proof?\\[-3ex]
 \end{enumerate}
 
A blend of mathematical reasoning supported by some silicon or quantum computation or a classical proof of excessive length and complexity (for example, the classification of finite simple groups) are examples of ``unconventional'' proofs. The ultimate goal of the mathematical activity is  the {\it advance human understanding of mathematics} (whatever this means!).
 
The answer to the first two question is affirmative. Indeed,  computations are represented in the  programming language used by the   computer
(the  `unconventional' computer too), even if the whole proof cannot be  globally `visualized' by a human being. The proof can be checked by any other mathematician having the equipment used in the 'unconventional' proof. 
A  proof   provides  knowledge only to the extent that its syntactic dimension is balanced by the semantic one; any gap
between them  makes the proof devoid of knowledge and paves the way for
 the proof to become a ritual without 
meaning. Proofs generating  knowledge, quite often  produce much more, for example,    'insight' (think of
the insight provided by understanding the algorithm used in the proof). 
  
 A misleading analogy would be  to replace,  in the above questions,  {\it `unconventional' methods}  with {\it ``testimony from a respected
 and (relevantly) competent mathematician''}. Certainly, such testimony  provides knowledge; it  does not qualify as a mathematical proof (even less   as a formalized proof), but the result is a ``mathematical activity'' because it advances our knowledge of mathematics.
 The  difference between `unconventional' methods and `relevant testimony'
can be found in the mechanisms used to produce outputs: a `relevant testimony' is the gut feeling of a respected, relevant, competent mathematician, by and large based on a considerable mathematical experience, while an `unconventional' method produces an objective argument.

There is little `intrinsic' difference between traditional and `unconventional' types of proofs as
 i) first and foremost,  {\it mathematical truth} cannot always be certified by proof, ii) correctness is not absolute, but almost certain, as mathematics advances by making mistakes and correcting and re--correcting them (mathematics fallibility was argued by Lakatos), iii) 
 non--deterministic and probabilistic
 proofs do not allow mistakes in the applications of rules, they are just indirect forms of checking (see Pollack \cite{pollack}, p. 210) which correspond to various degrees of rigour, iv) the  explanatory component, the understanding `emerging'  from proofs, while
 extremely important from a cognitive point of view, is subjective and has no bearing on formal correctness. As Hersh  noticed, mathematics like music exists by some logical, physical and biological manifestation, but ``it makes sense only as a mental and a cultural activity" (\cite{hersh}, p. 22).

How do we continue to
produce rigorous mathematics when more research will be performed in large computational
environments where we might or might not be
able to determine what the system has done or
why\footnote{Metaphorically described as ``relying on proof
by `Von Neumann says'".} is an open question.
The blend of logical and empirical--experimental
arguments are here to stay and develop. Of course,  some will continue to reject this trend, but, we believe,
 they will have as much effect as  King Canute's royal order to the tide. There are many reasons which support this prediction. They range from
economical ones (powerful computers will be more and more accessible to more and more people), social ones (skeptical oldsters are replaced naturally by youngsters born with the new technology, results and success inspire emulation) to pure mathematical (new challenging problems, wider perspective) and philosophical ones (note that  incompleteness is based on the analysis of the computer's behaviour).
The picture holds marvelous promises and challenges; it does not eliminate the continued importance of extended personal interactions in training and research.
  
\section*{Acknowledgements}
This paper is based on a talk presented at
the Workshop {\it Truths and Proofs},  a satellite meeting of the
{\it Annual Conference of the Australasian Association of Philosophy (New Zealand Division)}, 
Auckland, New
Zealand, 
December 2001. We are most grateful to Andreea Calude, Greg Chaitin, Sergiu Rudeanu, Karl Svozil, Garry Tee,  and Moshe Vardi for inspired comments  and suggestions.


\begin{thebibliography}{999}
\bibitem{ma} M. Aschbacher. The status of the classification of finite simple groups,
{\em Notices of the Amer. Math. Soc.}51, 7 (2004), 736--740.

\bibitem{auburn} D. Auburn. {\em Proof. A Play}, Faber and Faber, New York, 2001.

\bibitem{ah} K. Appel, W. Haken. {\em Every Planar Graph is Four Colorable},
Contemporary Mathematics 98, AMS, Providence,  1989.

\bibitem{Babai} L. Babai.  Probably true theorems, cry wolf? {\em Notices of the Amer. Math. Soc.} 41 (5) (1994), 453--454.

\bibitem{bb} D. H. Bailey, J. M. Borwein. Experimental mathematics: Recent developments
and future outlook, in B. Engquist, W. Schmid (eds.).  {\em World Mathematical
Year  Mathematics Unlimited---2001 and Beyond}, Springer-Verlag, Berlin, 2001,  51--66.


\bibitem{barrow} 
 J. Barrow. {\em Impossibility. The Limits of Science and the Science
of Limits}, Oxford University Press, Oxford, 1998.

\bibitem{bb1} J. M. Borwein, D. H. Bailey,    {\em The Experimental Mathematician.
Plausible Reasoning in the 21st Century,} 
A. K. Peters, Natick, Ma., 2003.


\bibitem{bbg} J. M. Borwein, D. H. Bailey,  R.  Girgensohn. {\em  Experimentation
in Mathematics. Computational Paths to Discovery}, A. K. Peters, Natick, Ma.,
2004.

\bibitem{beyer}  W. A. Beyer. The computer and
mathematics. {\em Notices of the Amer. Math. Soc.} 48, 11
(2001), 1302.

\bibitem{Blum} M. Blum. How to prove a theorem so no one else can claim it,
{\em Proceedings of the International Congress of Mathematicians,}
 Berkeley, California, USA, 1986, 1444--1451.

\bibitem{brent} R. P. Brent. Computation of the regular continued fraction for Euler's constant, {\em Math. of Computation} 31, 139 (1977), 771--777.


\bibitem{bridges} D. S. Bridges. Constructive truth in practice, in	 H.~G.~Dales and G.~Oliveri (eds.). {\em Truth in Mathematics},
Clarendon Press, Oxford, 1998, 53--69.

\bibitem{andreea} A. S. Calude. The journey of the four colour theorem through time,
{\em The NZ Mathematics Magazine} 38, 3 (2001), 27--35.
\bibitem{cris2002} C. S. Calude. {\em Information and Randomness: An Algorithmic
Perspective}, 
 2nd Edition,
Revised and Extended, Springer Verlag, Berlin,  2002.

\bibitem{crisomega} C.~S.~Calude. Chaitin $\Omega$ numbers, Solovay machines
and incompleteness,  {\em
Theoret. Comput. Sci.},  284 (2002), 269--277.


\bibitem{cris}  C.~S. Calude. Incompleteness, complexity, randomness and beyond,
{\em Minds and Machines: Journal for Artificial Intelligence, Philosophy and 
Cognitive Science}, 12, 4 (2002), 503--517.

\bibitem{cc}  C. S. Calude, G. J. Chaitin. Randomness everywhere, {\em Nature},
 400, 22 July (1999), 319--320.

\bibitem{crisds} C.~S. Calude, M.~J. Dinneen and C.-K. Shu.	Computing a glimpse of randomness,  {\em Experimental Mathematics} 11,
 2 (2002), 369--378.
 
 \bibitem{cj} C. S. Calude, H. J{\"u}rgensen. {\em Is Complexity a Source of Incompleteness?}, {\em CDMTCS Research Report} 241, 2004, 15 pp.
 
\bibitem{cjz} C.~Calude, H.~J\"{u}rgensen, M.~Zimand. Is independence an
exception\,?, {\em Appl. Math. Comput.} 66 (1994), 63--76.

\bibitem{crismarcus} C.~S. Calude, S. Marcus. Mathematical proofs at a crossroad? in
J. Karhum\" aki, H. Maurer, G. P\u aun, G. Rozenberg (eds.).
{\em Theory Is Forever},
%Essays Dedicated to Arto Salomaa on the Occasion of His 70th Birthday
Lectures
Notes in Comput. Sci. 3113, Springer Verlag, Berlin, 2004, 15--28.

%\bibitem{coins} C.~S. Calude, B. Pavlov. Coins, Quantum Measurements, and  Turing's %Barrier, {\em Quantum Information Processing}, 1, 1--2 (2002),  107--127. 

\bibitem{cp} C. S. Calude, G. P\u{a}un. {\em Computing with Cells and Atoms},
Taylor \&
Francis Publishers, London, 2001.


\bibitem{casti} J. L. Casti. {\em Mathematical Mountaintops}, Oxford University Press,
Oxford, 2001.


\bibitem{ch75} G. J. Chaitin. Randomness and mathematical proof,
{\em Scientific  American}, 232 (5) (1975), 47--52.

%\bibitem{ch99} G. J. Chaitin. {\em The Unknowable}, Springer Verlag,
%Singapore, 1999.

\bibitem{ch00} G. J. Chaitin. {\em Exploring Randomness}, Springer Verlag,
London, 2001.

\bibitem{ch02} G. J. Chaitin. Computers, paradoxes and the foundations of mathematics,
{\em American Scientist}, 90 March--April (2002), 164--171.


\bibitem{gregpccris}  G. J. Chaitin. Personal communication to C. S. Calude, 5 March, 2002.

\bibitem{gregphil} G. J. Chaitin. On the intelligibility of the universe and the notions of simplicity, complexity and irreducibility, {\tt http://www.cs.auckland.ac.nz/CDMTCS/ chaitin/bonn.html}, September 2002.

\bibitem{chomsky} N. Chomsky. {\em Syntactic Structures}, Mouton, The Hague, 1957.


\bibitem{cls} A. Connes, A. Linchnerowicz, M. P. Sch\" utzenberger. {\em Triangle of Thoughts},
AMS, Providence, 2001.

\bibitem{dijkstra} E. W. Dijkstra. Real mathematicians don't prove,
{\em EWD1012},   University of
Texas at Austin, 1988, {\tt http://www.cs.utexas.edu/users/EWD/EWD1012.pdf}.

\bibitem{deutsch-97} D. Deutsch. {\em The Fabric of Reality}, Allen Lane, Penguin Press,
1997.



\bibitem{demillo} R. A. De Millo, R. J. Lipton, A. J. Perlis. Social processes and proofs of
theorems and programs, {\em Comm. ACM} 22, 5 (1979), 271--280.

%\bibitem{devlin} K. Devlin. Lost innocence, {\em Guardian} 26 April (2001), 
%{\tt www.guardian.co.uk/Archive/}  {\tt Article/0,4273,4157878,00.html}.

\bibitem{dieu} J.  Dieudonn\'{e}.  {\em Pour L'honneur
de l'Esprit Humain},  Gallimard, Paris, 1986.

\bibitem{ft} W. Feit, J. G. Thomson. Solvability of groups of odd order,
{\em Pacific J. Math. } 13 (1963), 775--1029.

\if01
\bibitem{fenstad} J. E. Fenstad. Is mathematical still the science of paper, pencils and proofs?,
Conference on {\em Electronic Communication and Research in Europe},
Darmstadt/Seeheim, 15--17 April 1998,
{\tt //academia.darmstadt.gmd.de/seeheim/thebook/finals/fenstad.html}.
\fi

\bibitem{GMR}  S. Goldwasser, S. Micali,  C. Rackoff. The knowledge complexity of 
interactive proof--systems, {\em SIAM Journal of Computing,}
  18(1) (1989), 186--208.
  
\bibitem{foias} C. Foia\c{s}. Personal communication to S. Marcus (about 20 years ago).

\bibitem{hahn} L. E. Hahn, B. Epstein. {\em Classical Complex Analysis}, Sudury,
Mass., Jones and Barlettt, 1996.


\bibitem{harris} M. Harris. Contexts of justification, {\em The Mathematical Intelligencer}
23, 1 (2001), 10--22.


\bibitem{hayashi} S. Hayashi, R. Sumitomo, K. Shii. Towards the animation of proofs--testing
proofs by examples, {\em Theoret. Comput. Sci.} 272 (2002), 177--195.
\bibitem{hunt} R.  Hunt. The philosophy of proof, 
{\tt http://plus.maths.org/issue10/  features/proof4/}.
\bibitem{hersh} R. Hersh. {\em What Is Mathematics, Really\/?}, Vintage, London, 1997.


\bibitem{kline} M. Kline. {\em Mathematical Thought from Ancient to Modern Times}, Oxford University Press, Oxford, 
Vol. 1, 1972.
\bibitem{knuth} D. E. Knuth. Theory and practice, {\em EATCS Bull.} 27 (1985), 14--21.

%\bibitem{manin} Yu. I. Manin. Mathematics as profession and vocation, in
%V.~Arnold, M.~Atyah, P.~Lax,  B.~Mazur (eds.). {\em Mathematics: Frontiers 
%and Perspectives},
%IMU, AMS, Providence, 2000.



\bibitem{mc} D. MacKenzie. Slaying the kraken.  The sociohistory of a mathematical proof,
{\em Social Studies of Science} 29, 2 (1999), 7--60.

\bibitem{manin1} Yu. I. Manin.  Truth, rigour, and common sense, in	 H.~G.~Dales and G.~Oliveri (eds.). {\em Truth in Mathematics}, 
Clarendon Press, Oxford, 1998,  147--159.

\bibitem{marcus} S. Marcus.
No system can be improved in all respects, in  G. Altmann and
W. A. Koch (eds.). {\em Systems.
New Paradigms for the Human Sciences},  Walter de Gruyter, Berlin, 1998,
143--164.


\bibitem{marcus1} S. Marcus.
Imprecision, between variety and uniformity: the conjugate
pairs, in  J. J. Jadacki and  W.
Strawinski (eds.). {\em The World of Signs}, Poznan Studies in the Philosophy
of Sciences and the Humanities  62  Rodopi, Amsterdam, 1998 59--72.


\bibitem{marcus2} S. Marcus.  Metaphor as dictatorship, 
in  J. Bernard,
J. Wallmannsberger and  G. Withalm (eds.). {\em World of Signs. World of Things},  Angewandte Semiotik 15,  OGS,  Wien, 1997, 87--108.

\bibitem{pollack} R. Pollack. How to believe a machine--checked proof, in G. Sambin and
J. M. Smith (eds.). {\em Twenty--five Years of Constructive Type Theory}, Clarendon Press, Oxford,
1998, 205--220.

\bibitem{renyi} A. R\' enyi.
   {\em  Dialogues on Mathematics},  Holden Day, San Francisco, 1967.
   
\bibitem{rsst}  N. Robertson, D. Sanders, P. Seymour, R. Thomas. A new proof of the four--colour
theorem, {\em Electronic Research Announcements of the AMS} 2,1 (1996), 17--25.

\bibitem{rosen} R. R. Rosen.
Complementarity in social structures, {\em Journal of Social
and Biological Structures} 1 (1978), 191--200.

\bibitem{thom1}
R. Thom. Modern
mathematics: does it exist\/?, in  A.~G.~Howson (ed.). {\em Developments in Mathematical
Education}, Cambridge University Press, 1973, 194--209.

\bibitem{thom2} R. Thom. Topologie et linguistique, 
in A.~Haefliger and R.~Nerasimham (eds.). {\em   Essays on Topology and Related
Topics. Memoires
D\'edi\'es \`a Georges de Rham}, Springer Verlag, New York, 1970, 226--248.

\bibitem{solovay2k}
R. M. Solovay.  A version of $\Omega$ for which $ZFC$ can not predict a
single bit, in C. S. Calude and  G. P\u{a}un (eds.). {\em Finite Versus
Infinite. Contributions to an Eternal Dilemma}, 
Springer Verlag, London, 2000, 323--334.


\bibitem{tymoczko} T. Tymoczko. The four--colour problem and its philosophical significance,
{\em J. Philosophy} 2,2 (1979), 57--83.

\bibitem{schmidt}  U. Schmidt. {\em \"Uberpr\"ufung des Beweises f\"ur den
Vierfarben Satz}, Diplomarbeit Technische Hochschule, Aachen,
1982.


\bibitem{swart} E. R. Swart. The
philosophical implications of the four--colour problem, {\em American
Mathematical Monthly} 87, 9 (1980), 697--702.
\end{thebibliography}
 \end{document}